\documentclass{amsart}
\usepackage{amsmath, amsthm, amssymb}

\newtheorem{thm}{Theorem}

\theoremstyle{remark}
\newtheorem*{rem}{Remark}

\newcommand{\ZZ}{\mathbb{Z}}
\newcommand{\NN}{\mathbb{N}}
\newcommand{\QQ}{\mathbb{Q}}
\newcommand{\RR}{\mathbb{R}}
\newcommand{\CC}{\mathbb{C}}
\newcommand{\PP}{\mathbb{P}}
\renewcommand{\H}{\mathbb{H}}

\newcommand{\eqdef}{\stackrel{\text{def}}{=}}

\title{Asymptotic expansions, $L$-values and a new Quantum Modular Form}
\author{Edgar Costa}
\address[Edgar Costa]{Courant Institute of Mathematical Sciences \\
New York University \\
251 Mercer Street \\
New York, N.Y. 10012-1185 \\
U.S.A}
\thanks{The first author was partially supported by FCT doctoral grant SFRH/BD/69914/2010.}
\email{edgarcosta@nyu.edu}

\author{Korneel Debaene}
\address[Korneel Debaene]{Vakgroep Wiskunde \\
Ghent University \\
Krijgslaan 281, S22 \\
9000 Gent \\
Belgium}
\thanks{The second author was partially supported by the Fonds Professor Frans Wuytack.}
\email{korneel.debaene@ugent.be}

\author{Jo\~{a}o Guerreiro}
\address[Jo\~{a}o Guerreiro]{Department of Mathematics \\
Columbia University \\
Rm 509, MC 4406 \\
2900 Broadway \\
New York, NY 10027 \\
U.S.A}
\thanks{The third author was partially supported by FCT doctoral grant SFRH/BD/68772/2010.}
\email{guerreiro@math.columbia.edu}

\begin{document}

\begin{abstract}
In 2010 Zagier introduced the notion of a quantum modular form.
One of his first examples was the "strange" function $F(q)$ of Kontsevich.
Here we produce a new example of a quantum modular form by making use of some of Ramanujan's mock theta functions.
Using these functions and their transformation behaviour, we also compute asymptotic expansions similar to expansions of $F(q)$.
\end{abstract}

\maketitle
\section{Introduction and statement of results}

In 2010, Zagier \cite{zagier2} sketched a definition of a quantum modular form. Unlike classical modular forms, which are defined on the upper half plane, a quantum modular form is a function only defined on a subset of $\PP^1(\QQ)$. This set can be seen as the border of the upper half plane under the action of the modular group.
A quantum modular form should behave well under the action of some subgroup of the modular group, and should comply with some analytical constraints. 
Recall the action of an element $\gamma = \left( \begin{array}{cc}a&b\\c&d\end{array}\right) \in \mbox{SL}_2(\ZZ)$ on a function $f$,
\begin{equation*}
  f|_{\gamma}(z) = (cz+d)^{-k}f\left(\frac{az+b}{cz+d}\right).
\end{equation*}
Instead of demanding that $f$ is invariant under this action (which would give a vacuous definition), Zagier demands the difference
\begin{equation*}
f(z)-f|_{\gamma}(z)=h_\gamma(z)
\end{equation*}
to be a function which extends to an analytical function on $\PP^1(\RR)$.

The functions that appear in this paper have an integral period, either 1 or 2. Thus it is equivalent to view them as defined on $\H$ (or a vertical half-strip of breadth 1 or 2), or the open punctured unit disc. We will discern these two domains by exclusively using the variable $z$ on $\H$ and $q$ on the unit disc, the relationship between the two being either $q=e^{2\pi iz}$ or $q=e^{\pi iz}$. The domain of quantum modular forms, namely $\PP^1(\QQ)$ viewed as the border of the upper half plane, corresponds to the roots of unity and the origin of the unit disc.

One of the first examples studied by Zagier is the function known as Kontsevich's strange function \cite{zagier}
\begin{equation*}
F(q) \eqdef \sum_{n=0}^\infty \left(q;q\right)_n,
\end{equation*}
where we employ the standard abbreviation
\begin{equation*}
  (a;b)_n \eqdef \left( 1-a \left) \left( 1 -a b \left) \left( 1-ab^2 \right) \right. \right. \right. \right. \cdots \left( 1-ab^{n-1} \right), \; n \in \NN_0.
\end{equation*}
The function is a priori only defined when $q$ is a root of unity, but can be extended to a function in the open unit disc which satisfies near modular properties when we put $q = e^{2 \pi i z}$ with $z \in \mathbb{H}$.
Zagier also proves that $F$ has nice asymptotical expansions as $q$ tends to roots of unity radially. Furthermore, he establishes asymptotics for $F\left(\zeta_k\right)$ as $k\rightarrow\infty$, where $\zeta_k = e^{\frac{2\pi i}{k}}$ is a primitive $k$th root of unity.
Other asymptotical expansions of this kind have been studied, e.g. in \cite{rhoades} and \cite{ono}.

Our initial object of study is the function
\begin{equation*}
G(q)=\sum_{n=0}^\infty (-1)^n \left(q;q^2\right)_n,
\end{equation*}
which is only defined at odd roots of unity. However, this function can be extended to the open unit disc by a straightforward manipulation,
\begin{align}
    G(q) &= 1 + \sum_{n=1} ^\infty (-1)^n \left(q;q^2\right)_{n-1} \left( 1 - q^{2n -1}\right) \notag \\
    &= 1 + \sum_{n=1} ^\infty (-1)^n \left(q;q^2\right)_{n-1} -  \sum_{n=1} ^\infty q^{2n-1} (-1)^{n} \left(q;q^2\right)_{n-1} \notag \\
    &= 1 - G(q) + \sum_{n=0}^\infty  q^{2n+1} (-1)^{n} \left(q;q^2\right)_{n} \notag \\
    &= \frac{1}{2}\left( 1 +  \sum_{n=0} ^\infty q^{2n+1} \left( -1 \right)^n \left( q;q^2 \right) _n  \right) = \frac{1}{2} \phi(q), \label{1stform}
\end{align}
where $\phi(q)$ is one of Ramanujan's mock theta functions. Ramanujan defined it in the so-called Eulerian form as 
\begin{equation}
  \label{2ndform}
  \phi(q) \eqdef \sum_{n=0} ^\infty \frac{q^{n^2}}{\left(-q^2;q^2\right)_n}.
\end{equation}
The proof of the equality of the two forms of $\phi(q)$, (\ref{1stform}) and (\ref{2ndform}), is essential for the proof of \cite[Theorem 1.3]{ono}, where $\phi(q)$ appears intimately linked to $\psi(q)$, another mock theta function, 
\begin{equation*}
\psi(q) \eqdef \sum_{n=1} ^\infty \frac{q^{n^2}}{\left(q;q^2\right)_n} = \sum_{n \geq 0} q^{n+1} \left(-q^2;q^2\right)_n.
\end{equation*}

In this paper, we shall investigate how the properties of $G(q)$ are akin to those of $F(q)$, and of quantum modular forms, in general. Since the subject of quantum modular forms is still in its infancy, it is important to establish examples of such functions and determine their analytic properties.

In this paper, we prove the following statement.

\begin{thm}
The pair $\left(q^{-\frac{1}{24}}\phi(q),q^{-\frac{1}{24}}\psi(q)\right)$ is a vector-valued quantum modular form. More precisely,

\begin{enumerate}
    \item[1.] There are asymptotic expansions for $\phi\left(e^{-t} \zeta_{2k+1} ^l\right)$ and $\psi\left(e^{-t} \zeta_{4k} ^l\right)$ as $t \rightarrow 0^+$ (where $k,l \in \NN$ satisfying $(l,2k+1)=1$ and $(l,4k)=1$, respectively). In particular,
  \begin{equation*}
    \phi\left(e^{-t}\right) \sim \sum_{n=0}^\infty \frac{a_n}{n!} t^n \quad \mbox{as } t\rightarrow 0^+ ,
  \end{equation*}
  where
  \begin{equation*}
    a_n = \sum_{a + 2b + c = n} \frac{n!}{a! (2b)! c!} \left( \frac{3}{2} \right)^a \left( \frac{5}{2} \right)^{2b} E_{2a+2b} +  \sum_{a + 2b = n} \frac{n!}{a! (2b)!} \left( \frac{3}{2} \right)^a \left( \frac{1}{2}  \right)^{2b} E_{2a + 2b} ,
  \end{equation*}
$E_n$ are the Euler numbers, and the summations are taken over $a,b,c \in \NN_0$.

\item[2.] The functions $\phi(q)$ and $\psi(q)$ satisfy the following modular transformation equations:
\begin{equation*}  
  \begin{aligned}
  q^{- \frac{1}{24}} \phi(q) &= \sqrt{ \frac{4 \pi}{\alpha} } q_1 ^{- \frac{1}{24} } \psi\left(q_1\right) + \sqrt{ \frac{2 \pi}{3\alpha}} \int_0 ^\infty e^{ -\frac{\pi^2 u^2}{6 \alpha} } \frac{ \cosh \frac{5 \pi u}{6} + \cosh \frac{\pi u}{6}  }{\cosh \pi u}\, d u ,\quad \alpha = -\frac{ 2\pi i l}{2 k +1}, \\
  q^{- \frac{1}{24}} \psi(q) &= \sqrt{ \frac{\pi}{4 \alpha} } q_1 ^{- \frac{1}{24} } \phi\left(q_1\right) - \sqrt{ \frac{\pi}{6\alpha}} \int_0 ^\infty e^{ -\frac{\pi^2 u^2}{6 \alpha} } \frac{ \cosh \frac{5 \pi u}{6} + \cosh \frac{\pi u}{6}  }{\cosh \pi u}\, d u ,\quad \alpha = -\frac{2\pi i l}{ 4 k},
\end{aligned}
\end{equation*}
where $q = e^{-\alpha}$, $q_1 = e^{-\frac{\pi^2}{\alpha}}$, and $k$ and $l$ as above.

\item[3.] There are asymptotic expansions for $\phi\left(\zeta_{2k+1}\right)$ and $\psi\left(\zeta_{4k}\right)$ as $\NN \ni k \rightarrow \infty$,
\begin{align*}
    \zeta_{24(2k+1)} ^{-1} \phi\left(\zeta_{2k+1}\right) 
    &\sim \sqrt{2 (2k+1) i} \, \zeta_{96} ^{-23(2k + 1)} + \sum_{n=0} ^\infty \frac{\zeta_{24(2k+1)} ^{-1} b_n + \zeta_{24(2k+1)} ^{-25} c_n}{(2k+1)^n},\\
    \zeta_{96k} ^{-1} \psi\left(\zeta_{4k}\right) 
    &\sim  \sqrt{2ki} \, \zeta_{24}^{k} - \frac{1}{2} \sum_{n=0}^\infty \frac{\zeta_{96k} ^{-1} b_n + \zeta_{96k} ^{-25} c_n}{(4k)^n},
\end{align*}
where

\begin{align*}
  b_n &= \pi^n \sum_{a+2b = n} \frac{ (-1)^{a+b}(3i)^a}{a!(2b)!} E_{2a+2b}, & c_n &= \pi^n \sum_{a+2b = n} \frac{ (-1)^{a+b} (3i)^a 5^{2b} }{a!(2b)!} E_{2a+2b} ,
\end{align*}
 $E_n$ are the Euler numbers, and the summations are taken over $a,b \in \NN_0$.

\end{enumerate}

\end{thm}

\begin{rem}
  One remark on the modularity assertion is needed. As is usual in the theory of mock theta functions (see \cite{mcintosh}), we use the notation $q=e^{-\alpha}$, where $\Re \alpha > 0$. 
  Writing $\alpha = -\pi i z$, the transformation $q\mapsto q_1$ corresponds to the modular transformation $z\mapsto\frac{-1}{z}$. 
  Part 2 of the theorem then says that the vector-valued function  $f(z)=\left(q^{-\frac{1}{24}}\phi(q),q^{-\frac{1}{24}}\psi(q)\right)$ has the property that the vector \[f(z) - f\left(\frac{-1}{z}\right)\left(\begin{matrix} 
0 & \sqrt{\frac{i}{4 z}} \\
\sqrt{\frac{4 i}{z}} & 0 
\end{matrix}\right),\] 
  is given in terms of a certain improper integral, and extends to an analytic function on $\PP^1(\RR)$. Thus the pair $\left(q^{-\frac{1}{24}}\phi(q), q^{-\frac{1}{24}}\psi(q) \right)$ is a vector-valued quantum modular form with respect to the group
  \begin{equation*}
    \Gamma_\theta = \left\{ \begin{pmatrix}a & b \\c & d\end{pmatrix} \in \Gamma : b \equiv c \! \pmod{2} \right\},
  \end{equation*}
  which is generated by $z\mapsto\frac{-1}{z}$ and $z\mapsto z+2$, the latter being the trivial transformation.
\end{rem}

\section*{Acknowledgements}
We would like to thank professors Ken Ono and Rob C. Rhoades for bringing this problem to our attention. We would also like to thank the Arizona Winter School for creating opportunities for research and providing an excellent platform for collaboration.

\section{Nuts and bolts}
In this section we will gather the ingredients needed in the proof of our theorem. 

\subsection{Transformation formulae}

The near modular transformation behaviour of $\phi$ and $\psi$ on the open unit disc --- or equivalently the upper half plane via $q = e^{\pi i z}$, $z\in \mathbb{H}$ --- was found by Ramanujan and proved by Watson \cite{watson}. The following holds:
\begin{subequations}
\begin{align}
  q^{- \frac{1}{24}} \phi(q) &= \sqrt{ \frac{4 \pi}{\alpha} } q_1 ^{- \frac{1}{24} } \psi\left(q_1\right) + \sqrt{ \frac{6 \alpha}{\pi}} W(\alpha), \label{modulartransformationmcintosh1}\\
  q^{- \frac{1}{24}} \psi(q) &= \sqrt{ \frac{\pi}{4 \alpha} } q_1 ^{- \frac{1}{24} } \phi\left(q_1\right) - \sqrt{ \frac{3 \alpha}{2\pi}} W(\alpha), \label{modulartransformationmcintosh2}
\end{align}
\end{subequations}
where $q=e^{-\alpha}$, $\Re \alpha>0$, and $q_1 = e^{-\frac{\pi^2}{\alpha}}$. $W(\alpha)$ is defined as the following integral
\begin{equation}
  \label{W}
  W(\alpha) \eqdef \int_0 ^\infty e^{ -\frac{3}{2} \alpha x^2 } \frac{ \cosh \frac{5}{2} \alpha x + \cosh \frac{1}{2} \alpha x }{\cosh 3 \alpha x}\, d x. 
\end{equation} 

Many more of Ramanujan's mock theta functions satisfy similar curious transformation formulae, many of which were proven by Watson. A review can be found in \cite{mcintosh}.

One technical point is that the above expression for $W(\alpha)$ does not need to converge for $\Re \alpha = 0$. For $\Re \alpha >0$ changing variables $3 \alpha x = \pi u$ and moving the line of integration back to the real axis gives a new expression for $W(\alpha)$
\begin{equation}
  \label{Wextended}
  W(\alpha) =  \frac{\pi}{3 \alpha} \int_0 ^\infty e^{ -\frac{\pi^2 u^2}{6 \alpha} } \frac{ \cosh \frac{5 \pi u}{6} + \cosh \frac{\pi u}{6}  }{\cosh \pi u}\, d u,
\end{equation}  
which clearly converges for all purely imaginary $\alpha$ as well, and extends $W(t)$ continuously.

\subsection{Mordell Integral}

One important aspect of the above transformation formulas for $\phi$ and $\psi$ is the Mordell-type integral $W(\alpha)$ that appears as the obstruction to modularity. Zwegers \cite{zwegers} established in his thesis the modular properties of the Mordell integral.
For $z \in \CC$ and $\tau \in \mathbb{H}$ let
\begin{equation*}
  \label{h}
  h(z;\tau) \eqdef \int_\mathbb{R} \frac{e^{\pi i \tau x^2 - 2 \pi z x}}{\cosh \pi x} \, dx.
\end{equation*}
Then Zwegers proves
\begin{equation}
  \label{mordellmodular}
  h\left( \frac{z}{\tau} ; - \frac{1}{\tau}  \right) = \sqrt{ - i \tau } e^{- \frac{\pi i z^2}{\tau} } h\left( z;\tau \right).
\end{equation}

We may express $W(\alpha)$ using this notation:
\begin{equation}
  \label{Wintermsofh}
  W(\alpha) = \frac{\pi}{6 \alpha} \left( h\left( -\frac{1}{12} ; \frac{\pi i}{6 \alpha} \right) + h\left( -\frac{5}{12}; \frac{\pi i}{6 \alpha} \right) \right).
\end{equation}

\begin{rem}
  We could also obtain a modularity statement for $W(\alpha)$ by composing equations (\ref{modulartransformationmcintosh1}) and (\ref{modulartransformationmcintosh2}).
  The result is that $W(\frac{\pi^2}{\alpha}) = \left( \frac{\alpha}{\pi}\right)^{3/2} W(\alpha)$. However, using Zwegers function is a more general method since the function $h$ underlies the modularity properties of a wide class of integrals.
\end{rem}

\subsection{Euler Numbers}

The Euler numbers are the numbers which appear in the Taylor series expansion of the function $\frac{1}{\cosh x}$.
They will naturally enter the proof below due to the following equality (see \cite{euler}, for example):
\begin{equation}
  \label{euler}
  E_{2n} = (-1)^n 2^{2n} \int_\mathbb{R} \frac{w^{2n}}{\cosh \pi w} \, dw. 
\end{equation}

These numbers also show up in \cite{kac}, where they compute an asymptotic expansion of the Mordell integral in connection with the so-called \em{Kac-Wakimoto characters}.

\section{Proof of Theorem 1.1}

\begin{proof}[Proof of Theorem 1.1]
An interesting property of our mock theta functions is that the value at an (appropriate) root of unity exists by virtue of the fact that expression (\ref{1stform}) degenerates into a finite sum.
This is the key argument for the existence of the asymptotic expansions for the radial limits, i.e. limits of the form $\zeta_k^l e^{-t}$ with $t\rightarrow 0^+$.

We start by applying it to the simplest case $\phi\left(e^{-t}\right)$ with $t \rightarrow 0^{+}$.
Here $\phi(1)$ being expressed as a finite sum translates to

\begin{equation*}
  \left( e^{-t};e^{-2t} \right)_n = \left(1-e^{-t}\right) \left( 1 -e^{-3t} \right) \dots \left( 1-e^{-(2n+1)t} \right) = O\left(t^n\right).
\end{equation*}
Hence, each $k$th term in the asymptotic expansion can be expressed as a finite sum by expanding its first $k$th terms into their Taylor series. 

For $q = e^{-t} \zeta_{2k+1} ^l$ we have $\left(q;q^2 \right)_n = O\left(t^{\left\lfloor \frac{n+k}{2k+1} \right\rfloor} \right)$.
Therefore, each coefficient in the asymptotic expansion of $\phi\left(e^{-t} \zeta_{2k+1} ^l \right)$ as $t\rightarrow 0^+$ can also be obtained as a finite sum.
The same line of reasoning works for $\psi\left( e^{-t} \zeta_{4k} ^l \right).$

Furthermore, we can obtain an explicit asymptotic expansion for $q = e^{-t}$. Applying the transformation formula (\ref{modulartransformationmcintosh1}), which is valid for all $t>0$, we get
\begin{align*}
  e^{ \frac{t}{24}} \phi\left(e^{-t}\right) &= \sqrt{ \frac{4 \pi}{t} } e^{ \frac{\pi^2}{ 24 t} } \psi\left(e^{-\frac{\pi^2}{t}}\right) + \sqrt{ \frac{6 t}{\pi}} W(t).
\end{align*}
The first term of the right hand side is $O\left( e^{-\frac{O(1)}{t}} \right)$ and the second term has been computed (up to a factor of 2) in \cite[Proof of Theorem 4.1]{ono}, which settles the first part of the theorem.

The second part of the theorem is a question of taking limits of (\ref{modulartransformationmcintosh1}) and (\ref{modulartransformationmcintosh2}) to appropriate roots of units. We recall that the integral $W(\alpha)$ has been continuously extended to $\alpha \in \QQ i$ in Section 2.1. Note first how the roots of unity transform under the modular transformation $q\mapsto q_1$:
\begin{equation*}
\begin{aligned}
q = \zeta_{2k+1} ^l &\mapsto q_1 = \zeta_{4l} ^{-2k-1},\\
q = \zeta_{4k} ^l &\mapsto q_1 = \zeta_l ^{-k} \quad \quad (l \textup{  is odd}).
\end{aligned}
\end{equation*} 

This implies that the components of the transformation formulas (\ref{modulartransformationmcintosh1}) and (\ref{modulartransformationmcintosh2}) have well defined terms for $q$ an odd root of unity in (\ref{modulartransformationmcintosh1}) and for $q$ a $(4k)$th root of unity in (\ref{modulartransformationmcintosh2}), respectively.
Moreover, each of these terms is obtained by taking a limit; the Mordell integral has been continuously extended, and the other terms possess asymptotic expansions.
Consequently the validity of the equations (\ref{modulartransformationmcintosh1}) and (\ref{modulartransformationmcintosh2}) is kept while taking the limit to the desired roots of unity.

For the third part of the theorem, we apply the transformation formulas (\ref{modulartransformationmcintosh1}) and (\ref{modulartransformationmcintosh2}) and plug in $\zeta_{2k+1}$ and $\zeta_{4k}$, that is, we set $\alpha = -\frac{2\pi i}{2k+1}$ and $\alpha = -\frac{2\pi i}{4k}$, respectively:
\begin{equation*}
\begin{aligned}
  \zeta_{2k+1}^{-\frac{1}{24}} \phi\left(\zeta_{2k+1}\right) &= \sqrt{2(2k+1)i} \, \zeta_{4} ^{\frac{2k+1}{24} } \psi\left(\zeta_{4} ^{-2k-1}\right) + \sqrt{ \frac{-12 i }{2k+1}} W\left(-\frac{ 2\pi i }{2 k +1} \right), \\
  \zeta_{4k}^{- \frac{1}{24}} \psi\left(\zeta_{4k}\right) &= \sqrt{ \frac{k i}{2} } \zeta_{1}^{\frac{k}{24}} \phi(1) - \sqrt{ \frac{-3i}{4k}} W\left(-\frac{ 2\pi i }{4 k} \right).
\end{aligned}
\end{equation*}
We use that $\psi(\pm i) = \pm i$ and $\phi(1) = 2$, and get
\begin{equation*}
\begin{aligned}
\zeta_{24(2k+1)}^{-1} \phi\left(\zeta_{2k+1}\right) &= \sqrt{2(2k+1)i} \, \zeta_{96} ^{-23(2k+1)} + \sqrt{ \frac{-12 i}{2 k + 1}} W\left(-\frac{ 2\pi i }{2 k +1} \right), \\
  \zeta_{96k}^{-1} \psi\left(\zeta_{4k}\right) &= \sqrt{2ki} \, \zeta_{24}^{k} - \sqrt{ \frac{-3i}{4k}} W\left(-\frac{ 2\pi i }{4 k} \right).
\end{aligned}
\end{equation*}
The last ingredient we need is the asymptotics of $W\left(-\frac{2 \pi i}{m} \right)$. We use (\ref{Wintermsofh}) to express this in terms of $h(z; \tau)$:

\begin{equation*}
\begin{aligned}
  W\left( -\frac{2 \pi i}{m} \right) &= \frac{m i}{12} \left( h\left( - \frac{1}{12};-\frac{m}{12} \right) +  h\left( - \frac{5}{12};-\frac{m}{12} \right) \right) \\
  &= \sqrt{\frac{m i}{12}} \left( e^{-\frac{ \pi i}{12 m}} h\left( -\frac{1}{m} ; \frac{12}{m} \right) + e^{-\frac{ 25 \pi i}{12 m}} h\left( -\frac{5}{m} ; \frac{12}{m} \right) \right),
\end{aligned}
\end{equation*}
where in the last equation we applied (\ref{mordellmodular}) with $z=\frac{-1}{m}$ or $z=\frac{-5}{m}$ and $\tau=\frac{12}{m}$.

For $A\in\{-1,-5\}$,
\[h\left(\frac{A}{m};\frac{12}{m} \right)=\int_{\RR} \frac{e^{\pi i \frac{12}{m} x^2 - 2 \pi \frac{A}{m} x}}{\cosh(\pi x)} dx.\]
Expanding the term $e^{\pi i \frac{12}{m} x^2 - 2 \pi \frac{A}{m} x}$ in its Taylor series and applying identity (\ref{euler}), we get

\begin{equation*}
h\left(\frac{A}{m};\frac{12}{m}\right) \sim \sum_{n=0} ^{\infty} \frac{a(A)_n}{m^n},
\end{equation*}
where
\begin{equation*}
a(A)_n = \pi^n \sum_{a+2b = n} \frac{(-1)^{a+b} (3i)^a A^{2b} }{a!(2b)!} E_{2a+2b} .
\end{equation*}

Plugging in $A=-1$ and $A=-5$, we get the asymptotics of $W\left( -\frac{2 \pi i}{m} \right)$, for $m=2k+1$ or $m=4k$ and $k \rightarrow \infty$, thus finishing the proof of the theorem.
\end{proof}

\bibliographystyle{alpha}
\bibliography{CostaDebaeneGuerreiro}

\end{document}